\documentstyle[remark,12pt]{amsart}
\newbox\smilebox
\newbox\anchorbox
\newbox\noanchorbox
\newbox\tempbox

\setbox\smilebox=\hbox{$\smile$}

\def\anchor{\hbox{\vtop{
           \hbox to \wd\smilebox{\hfil\vrule width.4pt height7pt depth1pt\hfil}
           \vskip  -11.5truept
           \hbox to \wd\smilebox{\hfil$\smile$\hfil}}}}
\setbox\anchorbox=\anchor
\def\noanchor{\hbox{\vtop{
           \hbox to \wd\anchorbox{\hfil\anchor\hfil}
           \vskip -14truept
           \hbox to \wd\anchorbox{\hfil/\hfil}}}}
\setbox\noanchorbox=\noanchor

\def\fg#1#2#3{\setbox\tempbox=\hbox{$\scriptstyle{#2}$}
\ifnum\wd\anchorbox>\wd\tempbox\dimen255=\wd\anchorbox
\else\dimen255=\wd\tempbox\fi
{#1\,\vtop{\hbox to \dimen255{\hfil\anchor\hfil}
           \vskip -6truept
           \hbox to \dimen255{\hfil$\scriptstyle{#2}$\hfil}}
           \,#3}}

\def\nfg#1#2#3{\setbox\tempbox=\hbox{$\scriptstyle{#2}$}
\ifnum\wd\noanchorbox>\wd\tempbox\dimen255=\wd\noanchorbox
\else\dimen255=\wd\tempbox\fi
{#1\,\vtop{\hbox to \dimen255{\hfil\noanchor\hfil}
           \vskip -6truept
           \hbox to \dimen255{\hfil$\scriptstyle{#2}$\hfil}}
           \,#3}}

\setbox1=\hbox{$\bot$}

\def\north#1#2{#1\,
\hbox{$\bot$\llap {\hbox to\wd1 {\hfil $/$\hfil}}}
\,#2}

\def\nao#1#2#3{#1\  \hbox{\vtop{ 
\baselineskip=4pt
\hbox{$\bot$\llap {\hbox to\wd1 {\hfil $/$\hfil}}
\hskip .05em \llap{\hbox{$^{\scriptscriptstyle{a}}$}}}\hbox{$\scriptstyle
{#2}$}}}\, #3}

\def\bp{\par{\bf Proof.}$\ \ $}

\def\includeE#1{{\lhook\kern-3.5pt\joinrel\smash{
    \mathop{\longrightarrow}\limits^{#1}}}}

\def\efor/{Example~\ref{E4}}

\def\ep{\par\bigskip}

\def\BL/{Baldwin--Lachlan}
\def\Bu/{Buechler}
\def\Hr/{Hrushovski}
\def\lm/{locally modular}
\def\wm/{weakly minimal}
\def\nm/{non--modular}
\def\tt/{totally transcendental}
\def\ud/{unidimensional}
\def\sm/{strongly minimal}

\def\abar{\overline{a}}

\def\bbar{\overline{b}}

\def\hbar{\overline{h}}

\def\xbar{\overline{x}}
\def\ybar{\overline{y}}

\def\delbar{\overline{\delta}}

\def\tao{\tau}

\def\dom{{\rm dom}}

\def\tr/{trivial}
\def\nt/{non--trivial}
\def\st/{strong type}

\def\TV/{Tarski--Vaught}

\def\sc/{sound construction}
\def\ac/{atomic construction}

\def\fal/{functional}
\def\upl/{unique parallel lines}
\def\chp/{categorical in a higher power}

\def\text#1{\ \hbox{#1}\ }

\def\sbq{\subset}
\def\contains{\supseteq}
\def\<{\langle}
\def\>{\rangle}
\def\conc{\widehat{~~}}
\def\forces{\mathrel{\raise.4ex\hbox{$\scriptstyle \vert$\hskip-.5ex}\vdash}}

\def\K/{${\cal K}$}
\def\KM{{\cal K}}

\def\r{\restriction}
\def\rl{\r L}
\def\KMOD/{$\KM={\rm Mod}(T_1)\rl$}
\def\LK/{$L(\KM)$}
\def\S{{\rm Spec}(\KM)}
\def\al{\alpha}
\def\b{\beta}
\def\d{\delta}

\def\e{\emptyset}
\def\g{\gamma}
\def\k{\kappa}
\def\l{\lambda}
\def\o{\omega}

\def\SN/{\hbox{$\S\ne\e$}}
\def\ao{\aleph_0}
\def\a1{\aleph_1}

\def\PCD/{${\rm PC}_{\a0}$}
\def\PC/{${\rm PC}_{\Delta}$}
\def\ra{\rightarrow}

\def\bg{\beth_\g}

\def\SG/{\hbox{$\S\cap\bg$}}
\def\SNG/{\hbox{$\S\cap\bg\not =\e$}}

\def\T{{\rm KT}}

\def\phi{\varphi}
\def\Ibar{\overline{I}}
\def\Jbar{\overline{J}}
\def\lan{\langle}
\def\ran{\rangle}
\def\cof{\hbox{cof}}
\def\coi{\hbox{coi}}
\def\endproof{\enspace\vrule height6pt width4pt depth0pt\ep} 
\def\sbq{\subseteq}

\def\Khom{{\bf K}_{\hbox{2tr}}}
\def\FF{{\cal F}}

\def\KK{{\bf K}}
\def\QQ{{\Bbb Q}}

\def\S{{\cal S}}

\def\f{\tilde}
\def\ff{\f{f}}
\def\range{{\rm range}}

\def\T{{\cal T}}
\newtheorem{theorem}{Theorem}[section]

\newtheorem{lemma}[theorem]{Lemma}
\newtheorem{prop}[theorem]{Proposition}

\newremark{defn}[theorem]{Definition}
\newremark{example}[theorem]{Example}
\newremark{remark}[theorem]{Remark}
\author{M. C. Laskowski}
\thanks{\hbox{Partially supported by  NSF 
Research~Grants~DMS~9403701 and DMS~9704364}}
\address{Department of Mathematics 
\\ University of Maryland \\ College Park, MD 20742}
\email{mcl@@math.umd.edu}
\author{S. Shelah}
\address{Department of Mathematics\\
Hebrew University of Jerusalem
\and Department of Mathematics\\
Rutgers University}
\thanks{The authors thank the U.S.-Israel
Binational Science Foundation for its support of this project.
This is item 560 in Shelah's bibliography.}
\date{\today}
\title{The Karp complexity of unstable classes}


\begin{document}

\begin{abstract}
A class $\KK$ of structures is {\it controlled\/} if, for all cardinals $\l$,
the relation of $L_{\infty,\l}$-equivalence partitions $\KK$ into a set of equivalence
classes (as opposed to a proper class).  We prove that the class of doubly transitive
linear orders is controlled, while any pseudo-elementary class with the $\o$-independence
property is not controlled.
\end{abstract}

\maketitle
%\newpage

\section{Introduction}
One of the major accomplishments of model theory has been the 
discovery of a dividing
line between those theories in a countable language whose models can be described
up to isomorphism by a reasonable set of invariants and those whose models cannot be so described.
Models of {\it classifiable\/} theories are described up to isomorphism by an
`independent tree' of countable elementary submodels, while the isomorphism
type of any unclassifiable theory cannot be described by any reasonable set of invariants
(see~\cite{Sh:c}).
Unfortunately, the great majority of
classes of structures studied
in mathematics are unstable, and thus
fall on the `non-structure' side of this divide.  Thus, it is desirable
to search for dividing lines between unstable classes of structures.
Our thesis is that while an 
unstable (pseudo-elementary)
class necessarily has the maximal number of non-isomorphic
models in every uncountable cardinality, 
it is still possible to assign a {\it set of invariants\/} to
some unstable classes of structures.  
In some cases (see e.g., Example~\ref{many})
the large number of non-isomorphic models is due simply to our
ability to code arbitrary stationary sets into the skeletons of
Ehrenfeucht-Mostowski models.
In other words, for some classes of structures
the reason for the non-isomorphism of two structures in the class
need not be very robust. Indeed, in such cases the structures can be forced
to be isomorphic by a forcing that
merely adds a new closed, unbounded subset of some cardinal
to the universe.
That is, although they are nonisomorphic,
the structures are not very different from each other.
On the other hand, for other classes of structures
(see Theorem~\ref{controlled}) there
are more serious obstructions to a structure theorem.

Our ultimate goal is to determine to which unstable classes of structures one can
associate a reasonable set of structural invariants.  These 
invariants need not (and typically will not) determine the structures up to isomorphism.
Instead, we ask that any two structures with the same invariants be very much the same.
In this paper we focus on  $L_{\infty,\l}$-equivalence
for various cardinals $\l$ and ask which unstable classes are partitioned into only a {\it set\/}
of equivalence classes (as opposed to a proper class).
We call a class $\KK$ {\it controlled\/} if $\KK$ has only a set of $L_{\infty,\l}$-equivalence
classes for all cardinals $\l$.
Typically, $L_{\infty,\l}$-equivalence does not characterize models up to isomorphism
even when we fix the cardinalities of the models.
(In \cite{Sh:220} the 
second author shows that for any unstable pseudo-elementary class
and any uncountable regular cardinal $\l$,
there are $2^\l$ non-isomorphic models of size $\l$ that are $L_{\infty,\l}$-equivalent.)
However, in some sense
two $L_{\infty,\l}$-equivalent structures of the same cardinality are very much the same.
For instance, if one uses the back-and-forth system 
witnessing their equivalence
as a notion of forcing,
then the two structures will become isomorphic in the corresponding
forcing extension.

In this paper we 
obtain two complementary results.  On one hand,
in Section~\ref{2transitive}
we analyze the pseudo-elementary class $\Khom$ of doubly transitive linear orders.
This class is unstable, hence the stigma of non-structure applies.
Despite this,
we prove that $KC_{\mu^+}(\Khom)\le\o$  (see Definition~\ref{class})
for all
uncountable cardinals $\mu$, hence $\Khom$ is controlled.
This is one of very few theorems in which an unstable pseudo-elementary
class shows any sign of structure.
On the other hand, in Section~\ref{omegaindep} we prove that any pseudo-elementary
class with the $\o$-independence property (see Definition~\ref{o-indep})
is not controlled.  In fact, if the language used in describing $\KK$ is
countable then
$KC_\l(\KK)=\infty$ for all cardinals $\l\ge\aleph_3$.

There is still much that we do not know about the notion of control.
A fundamental question that remains open is whether there is an unstable {\it elementary\/}
class that is controlled.
We conjecture, and hope to prove,
that any pseudo-elementary class with the independence property
is not controlled; this would substantially strengthen our second result.

\section{Controlled classes}

In this section we state a series of  definitions that lead to the concept of a class
of structures being
controlled (see Definition~\ref{controlled}).
We apply these definitions to the  theory of  dense linear orders to illustrate why
it is desirable to consider the $\l$-Karp complexity of a class for uncountable cardinals $\l$.
We first reintroduce the notion of a partial isomorphism, but with a slight
variation.
As we are only concerned with the definable subsets of structures (and not
their quantifier complexity) we insist that 
all partial isomorphisms are elementary maps.

\begin{defn}  Given two elementarily equivalent
structures $M$ and $N$ in the same language and an infinite cardinal $\l$,
a {\it $\l$-partial isomorphism\/} is a partial elementary map
with domain of cardinality less than $\lambda$, that is:
a function $f$ from a subset $D$ of $M$
into $N$
of size less than $\l$
satisfying
$$M\models\phi(d_1\dots d_n)\qquad \hbox{if and only if}\qquad N\models
\phi(f(d_1),\dots,f(d_n))$$
for all formulas $\phi(x_1,\dots,x_n)$ of the language and all $d_1,\dots,d_n$ from $D$.
We denote the family of $\l$-partial isomorphisms by $\FF_\l(M,N)$.  If 
$M=N$ we simply write $\FF_\l(M)$.
\end{defn}

The complexity
of $\FF_\l(M)$ is a measure of how deeply one needs to look to 
understand the relationship of a small 
subset (i.e., of size less than
$\l$) with the rest of the model.
In order to measure this depth we endow the family with the following rank.
\begin{defn}  \label{Rank}
For $f\in \FF_\l(M,N)$,
\begin{enumerate}
\item Rank$(f)\ge 0$ always;
\item For $\al$ limit, Rank$(f)\ge\al$ if and only if Rank$(f)\ge\b$ for all $\b<\al$;
\item Rank$(f)\ge\al + 1$ if and only if 
\begin{enumerate}
\item for all
$C\sbq M$ of size less than $\l$,
there is $g\in\FF_\l(M,N)$ extending $f$ with $C\sbq\dom(g)$
and Rank$(g)\ge\al$; and
\item dually, for all $C\sbq N$ of size less than $\l$,
there is $g\in\FF_\l(M,N)$ extending $f$ with $C\sbq\range(g)$
and Rank$(g)\ge\al$.
\end{enumerate}
\end{enumerate}
The {\it $\l$-Karp complexity\/} $KC_\l(M,N)$ of the pair of structures $M,N$ is the least
ordinal $\al$ such that Rank$(f)\ge\al$ implies Rank$(f)\ge\al+1$ for all $f\in\FF_\l(M,N)$.
Again, if $M=N$ we simply write $KC_\l(M)$.  
\end{defn}

The $\l$-Karp complexity 
of a structure is related to the notions of $L_{\infty,\l}$-Scott
height and back-and-forth systems.
It is a routine diagram-chasing  exercise to show that 
if two structures $M$ and $N$ are $L_{\infty,\lambda}$-equivalent (hence there is a back-and-forth
system in $\FF_\l(M,N)$) then $KC_\l(M)=KC_\l(N)$.

If one fixes the signature, then for any cardinals $\k$ and $\l$
it is easy to find an ordinal bounding the $\l$-Karp complexity of any structure 
of that signature of size at most $\k$.
By contrast, whether or not there is a upper bound on the $\l$-Karp complexities of 
all structures in a class $\KK$ that does not depend on $\k$
provides a robust dichotomy between classes.  
This is demonstrated by
the following definition and proposition.
The reader is referred to \cite{Barwise} for the undefined notions.

\begin{defn}  \label{class}
For $\KK$ a class of structures, the {\it $\l$-Karp complexity\/} of $\KK$, written $KC_\l(\KK)$,
is the supremum of the ordinals $KC_\l(M)$ among all  $M\in \KK$ if the supremum exists.
Otherwise, we set $KC_\l(\KK)=\infty$.
\end{defn}

%\begin{defn}  Two elementarily equivalent structures $M$ and $N$
%are {\it $(\l,\al)$-equivalent\/}, written $M\equiv_{\l,\al} N$, if
%the empty function in $\FF_\l(M,N)$ has Rank at least $\al$.
%\end{defn}

\begin{prop}  \label{equiv}
The following conditions are equivalent for a class of structures
$\KK$ and an infinite  cardinal $\l$.
\begin{enumerate}
\item $KC_\l(\KK)<\infty$;
\item The relation of $L_{\infty,\l}$-equivalence on $\KK$
has only a set of equivalence classes;
\item There are only a set of  $L_{\infty,\l}$-types of subsets of size
less than $\l$ realized in elements of $\KK$;
\item There are only a set of distinct $L_{\infty,\l}$-Scott sentences among the
elements of $\KK$;
\item There is a cardinal $\k$ such that the notions of 
$L_{\k,\l}$-equivalence and $L_{\infty,\l}$-equivalence coincide on $\KK$.
\end{enumerate}
\end{prop}

\bp  The implication $(2)\Rightarrow (1)$ follows from the observation that
$\l$-Karp complexity is preserved under $L_{\infty,\l}$-equivalence.  The implications
$(1)\Rightarrow(4)\Rightarrow(5)\Rightarrow(3)\Rightarrow(2)$ all follow easily.
\endproof

When  $\l=\ao$ the $\l$-Karp complexity often does not
yield much information about the inherent complexity of a class $\KK$.
For example, if $\KK$ is the class of models of an $\ao$-categorical theory, then 
$KC_{\ao}(\KK)=0$ since every model is $\ao$-homogeneous.
However, our thesis is that for larger $\l$, $\l$-Karp complexity gives a good measure
of the complexity of the class.
It follows from Proposition~\ref{equiv}(3) that if $KC_\l(\KK)=\infty$ for some
cardinal $\l$, then $KC_\k(\KK)=\infty$ for all larger cardinals $\k$.
This leads us to the crucial definition of the paper.

\begin{defn} \label{controlled}
A class  $\KK$ of structures is  {\it controlled\/}
if $KC_\l(\KK)<\infty$ for all infinite cardinals $\l$.
\end{defn}

Note that if a class $\KK$ is controlled, then it follows from
Proposition~\ref{equiv}(2) that
for every cardinal $\l$,
the relation of $L_{\infty,\l}$-equivalence partitions $\KK$ into only a set of
equivalence classes (as opposed to a proper class). 
Continuing our example, $KC_{\ao}(DLO)=0$, as
$DLO$, the theory of dense linear orders with no endpoints
is $\ao$-categorical.
However, this observation hides the
fact that one can code arbitrary ordinals into 
dense linear orders.  
This ability to code ordinals implies that
the class DLO is not controlled.
In fact, $KC_\l(DLO)=\infty$ for all uncountable cardinals $\l$.
To see this, fix an uncountable cardinal $\l$ and, for each non-zero ordinal $\alpha$,
let $J_\alpha$ be the linear order with universe $(\eta\cdot\l)\cdot\alpha$,
where $\eta$ denotes the order type of the rationals.
In light of 
Proposition~\ref{equiv}(2)
it suffices to show that 
$J_\alpha$ is not $L_{\infty,\l}$-equivalent to $J_\beta$ whenever
$\alpha\neq\beta$.  
So choose non-zero ordinals $\alpha$ and $\beta$
such that $J_\alpha$ is $L_{\infty,\l}$-equivalent to $J_\beta$. Let
$E$ be the equivalence relation such that
$E(x,y)$ if and only if there are fewer than $\l$ elements between $x$ and $y$.
Since $E$ is expressible in the logic $L_{\infty,\l}$, this implies that
the condensation $J_\alpha/E$ is $L_{\infty,\l}$-equivalent to $J_\beta/E$.
But $(J_\alpha/E,\le)\simeq (\alpha,\le)$, $(J_\beta/E,\le)\simeq(\beta,\le)$,
and it is readily checked that distinct ordinals are not even
$L_{\infty,\omega}$-equivalent.
Hence $\alpha$ must equal $\beta$.

\section{Doubly transitive linear orders}  \label{2transitive}

In this section we investigate the class $\Khom$ of
infinite doubly transitive linear orders.  That is, $(I,\le)\in\Khom$
if and only if the linear order $I$ is dense with no endpoints and for all
pairs $a<b$, $c<d$ from $I$, the interval $[a,b]$ is isomorphic to the interval $[c,d]$.
Such orders arise naturally:  The underlying
linear order of any ordered field
is necessarily doubly transitive.
Clearly, there is only one countable structure in $\Khom$
up to isomorphism.
The class $\Khom$ is a pseudo-elementary (PC) class
that is visibly unstable, so by \cite{Sh:c} there are $2^\l$
non-isomorphic structures in $\Khom$ of size $\l$ for all
uncountable cardinals $\l$.  Further, by \cite{Sh:220}, for all uncountable regular
cardinals $\l$
there is a family of $2^\l$ structures
in $\Khom$ of size $\l$ that are $L_{\infty,\l}$-equivalent, yet pairwise non-embeddable.

Nonetheless, the class of doubly transitive
linear orders is not entirely without structure.
There are natural `invariants' one
can associate with such orders.
These invariants will not determine the orders up to isomorphism,
but they will be sufficient to demonstrate that the 
$\l$-Karp complexity of $\Khom$ is bounded for all cardinals $\l$.

The most natural invariant of a doubly transitive linear order is the isomorphism
type of its closed intervals.  Accordingly, we call $I_0,I_1\in\Khom$
{\it locally isomorphic\/} and write $$I_0\sim I_1$$
if $[a,b]\simeq [c,d]$
for $a<b$ from $I_0$ and $c<d$ from $I_1$.
Evidently local isomorphism is an equivalence relation on $\Khom$ and
$I\sim J$ for any infinite convex subset $J\sbq I$, if $I\in\Khom$.

The second invariant was developed by Droste and Shelah in 
\cite{DrSh:195}.  The definitions that
follow are slight adaptations of similar notions used there.
The most notable variation is that in \cite{DrSh:195} there is no bound on the number of
levels of the {\it decomposition tree\/} and the cardinals $\l_\eta$ can
be any uncountable regular cardinal.

For the whole of this section,
fix an uncountable cardinal $\mu$.

\begin{defn} A {\it $\mu$-decomposition tree\/} is a 
subtree  $T$ of $\bigcup\{{^\al \mu}:\al<\mu^+\}$ satisfying:
\begin{enumerate}
\item $T$ is downward closed, i.e., $\eta\in T$ implies 
$\eta|\al\in T$ for all $\al<lg(\eta)$;

\item If $lg(\eta)$ is  a limit ordinal or 0 and $\eta|\al\in T$ for all
$\al<lg(\eta)$ then
$\eta\in T$ and $\eta$ has exactly two immediate successors; more specifically,
we require
$Succ_T(\eta)=\{\eta\conc\lan 0\ran,\eta\conc\lan 1\ran\}$;

\item If $\eta\in T$ and $lg(\eta)$ is a successor ordinal, then
either $Succ_T(\eta)=\emptyset$ or 
$Succ_T(\eta)=\{\eta\conc\lan\al\ran:\al\in C\}$
for some club subset $C$ of
a regular cardinal $\l_\eta\in [\aleph_1,\mu]$.
\end{enumerate}
Let $T^*=\{\eta\in T:lg(\eta)$ is a successor ordinal$\}$.
\end{defn}

We define a linear order on $T^*$ which is a cross between
lexicographic and antilexicographic order.
To every node $\eta$ of $T^*$ we first
associate a direction $dir(\eta)\in
\{{\rm LEFT,RIGHT}\}$.  Suppose $lg(\eta)=\d+n$,
where $\d$ is a limit ordinal or 0 and $n\in\o$.  Then
\begin{itemize}
\item $dir(\eta)={\rm LEFT}$ if $\eta(\d)+n$ is even;
\item $dir(\eta)={\rm RIGHT}$ if $\eta(\d)+n$ is odd.
\end{itemize}
The idea is that if $dir(\eta)=$LEFT, then the successors of $\eta$ will
all be to the left of $\eta$.  Each of these successors will have direction RIGHT,
so their successors will be to their right and so forth.  Formally,
the linear order $<^{T^*}$ is defined by the following clauses.
\begin{itemize}
\item  If $\eta\lhd \nu$ then $\eta <^{T^*}\nu$ if and only if $dir(\eta)=$RIGHT;
\item  If $\eta,\nu$ are incomparable, let $\g$ be least such that $\eta(\g)\neq\nu(\g)$
and let $\rho=\eta|\g$.
\begin{itemize}
\item  If $\g$ is a limit ordinal  or 0 
then $\eta <^{T^*}\nu$ if and only if $\eta(\g)=0$ and $\nu(\g)=1$;
\item  If $\g$ is a successor ordinal (so $\rho\in T^*$) and $dir(\rho)=$LEFT then $\eta <^{T^*}\nu$
if and only if $\eta(\g)<\nu(\g)$;
\item  If $\g$ is a successor ordinal and $dir(\rho)=$RIGHT, then $\eta <^{T^*}\nu$
if and only if $\eta(\g)>\nu(\g)$.
\end{itemize}
\end{itemize}

The following definition differs slightly from normal usage as we
include the endpoints.
\begin{defn}  For $I$ a dense linear order,  the {\it Dedekind completion\/}
of $I$ is the linear order $(\Ibar,\le^{\Ibar})$ 
with universe $$\Ibar=\{A\sbq I:A\ \hbox{downward closed with no largest element}\}$$
and $A\le^{\Ibar} B$ if and only if $A\sbq B$.
We let $-\infty$ denote the smallest element of $\Ibar$ and $+\infty$ denote the largest.
To simplify notation we identify the element $a\in I$ with $\{x\in I:x<a\}\in \Ibar$
and write e.g., $I\sbq \Ibar$.
If $J$ is a convex subset of $I$, then $\Jbar$ denotes the smallest closed interval
in $\Ibar$ that contains $J$
and we identify $\Jbar$ with the Dedekind completion of $J$.
\end{defn}

\begin{defn}  A {\it $\mu$-representation\/}
of a linear ordering $I$ is a pair $(T,g)$, where $T$ is a 
$\mu$-decomposition tree and $g:T^*\ra \Ibar$
is an order-preserving  function satisfying the following conditions:

\begin{description}
\item[1] $g(\lan 0\ran)=-\infty$, $g(\lan 1\ran)=+\infty$;
\item[2] If $lg(\eta)=\g+1$, where $\g>0$ is a limit ordinal, let
$D$ be the largest interval $[a,b]$ of $\Ibar$ such that for all successor
ordinals $\al<\g$, $D$ is between $g(\eta|\al)$ and $g(\eta|\al+1)$.
\begin{enumerate}
\item If $\eta(\g)=0$ then $g(\eta)=a$;
\item If $\eta(\g)=1$ then $g(\eta)=b$;
\item If $a=b$ then we call $\eta$ {\it degenerate}.
\end{enumerate}

\item[3] If $dir(\eta)$=LEFT then
\begin{enumerate}
\item $\eta$ is maximal in $T$ if and only if one of the three conditions
hold:
\begin{enumerate}
\item $\eta$ is degenerate;
\item $\cof(g(\eta))=\ao$;
\item $\cof(g(\eta))>\mu$;
\end{enumerate}
\item If $\eta$ is not maximal in $T$, then $Succ_T(\eta)=\{\eta\conc
\lan\al\ran:\al\in C\}$ for some club subset of $\cof(g(\eta))$,
and $\{g(\eta\conc\lan\al\ran):\al\in C\}$ is  continuous, strictly
increasing, and has supremum $g(\eta)$.
\end{enumerate}
\item[$3^*$] If $dir(\eta)=$RIGHT then
\begin{enumerate}
\item $\eta$ is maximal in $T$ if and only if one of the three conditions
hold:
\begin{enumerate}
\item $\eta$ is degenerate;
\item $\coi(g(\eta))=\ao$;
\item $\coi(g(\eta))>\mu$;
\end{enumerate}
\item If $\eta$ is not maximal in $T$, then $Succ_T(\eta)=\{\eta\conc
\lan\al\ran:\al\in C\}$ for some club subset of $\cof(g(\eta))$,
and $\{g(\eta\conc\lan\al\ran):\al\in C\}$ is  continuous, strictly
decreasing, and has infimum $g(\eta)$.
\end{enumerate}

\end{description}
\end{defn}

A $\mu$-representation $(T,g)$ partitions 
$\Ibar$ into a set of intervals $\{I_\eta:\eta\in T^*\}$
where
$I_{\eta\conc \lan 0\ran}=I_{\eta\conc \lan 1\ran}=
(g(\eta\conc \lan 0\ran),g(\eta\conc \lan 1\ran))$
for all $\eta\in T$ of limit length, and if 
$Succ_T(\eta)=\{\eta\conc\lan\al\ran:\al\in C\}$ for a club $C$
then 
$$I_{\eta\conc\lan\al\ran}=
\left\{ \begin{array}{ll}
(g(\eta\conc\lan\al\ran),g(\eta\conc\lan\al^+\ran)) &
                           \mbox{if $dir(\eta)=$LEFT;} \\
(g(\eta\conc\lan\al^+\ran),g(\eta\conc\lan\al\ran)) &
                           \mbox{if $dir(\eta)=$RIGHT}
       \end{array}\right. $$
where $\al^+$ is the least element of $C$ larger than $\al$.
It is easily shown by induction that the intervals $\{I_\eta:\eta\in T^*
\cap {^\al\mu}\}$ are pairwise disjoint for any fixed successor ordinal $\al$.

For any dense linear order $I$, one can build a $\mu$-representation $(T,g)$
of $I$ level by level by successively choosing a continuous, strictly increasing
[or decreasing]
sequence $\lan g(\eta\conc\lan\al\ran):\al\in \l_\eta\ran$ from the interval $I_\eta$.
At first blush, it appears that one has considerable freedom in such a construction.
However, our freedom is considerably
limited by the following observation.

\medskip\par\noindent
{\bf Observation}  Let $J$ be any linear order of cofinality $\l>\ao$.
For any club subsets $C_1,C_2$ of $\l$ and any two continuous, strictly
increasing, cofinal sequences  $\lan a_i:i\in C_1\ran$ and $\lan b_i:
i\in C_2\ran$ in $J$, the set $D=\{i\in C_1\cap C_2:a_i=b_i\}$ is a club
subset of $\l$.

\medskip
By repeatedly applying this observation to a pair of 
$\mu$-representations 
of a linear order, we see that they must  `agree on a club.' 
More precisely, call a subtree $T'$ of a $\mu$-decomposition tree $T$ a {\it club subtree\/}
if $T'$ itself is a $\mu$-decomposition tree and, for each
$\eta\in T'$ that is not maximal in $T'$, $Succ_T(\eta)$ and $Succ_{T'}(\eta)$
are both indexed by club subsets of the same regular cardinal.
If $(T_1,g_1)$ and $(T_2,g_2)$ are two $\mu$-representations of $I$,
then by using the observation above at each node
there is a $\mu$-representation $(T,g)$ of $I$ such that
$T$ is a club subtree of both $T_1$ and $T_2$ with 
$g(\eta)=g_1(\eta)=g_2(\eta)$
for all $\eta\in T^*$.
More generally we have the following definition and lemma.

\begin{defn}  A subset $A$ of a $\mu$-decomposition tree $T$ is {\it closed\/}
if
$A$ is downward closed, (i.e., if $\eta\in A$ then
$\eta|\al\in A$ for all $\al< lg(\eta)$) and
$A$ is closed under successor, (i.e., if $\eta\in A$ then
$Succ_T(\eta)\sbq A$). 
\end{defn}

Note that for any subset $A\sbq T$ of size at most $\mu$, there is a closed
subset $B\contains A$ of size at most $\mu$.

\begin{lemma}  \label{closed}
Suppose $(T,g)$ is a $\mu$-representation of $I_0$,
$S\sbq T$ is closed, and
$f_0,f_1:\Ibar_0\ra \Ibar_1$ are order-preserving, continuous partial
functions whose domains contain $\{g(\eta):\eta\in S\cap T^*\}$
that satisfy $f_0(-\infty)=f_1(-\infty)$ and $f_0(+\infty)=f_1(+\infty)$.
Then there is a club subtree $Y\sbq T$ such that
$$f_0(g(\eta))=f_1(g(\eta))$$
for all $\eta\in S\cap Y^*$.
\end{lemma}

\bp
We construct $Y$ by induction on the levels of $T$.
Assume that we have found $Y_\g$, a club subtree of $T\cap\bigcup\{{^\b\mu}:\b<\g\}$
such that $f_0(g(\eta))=f_1(g(\eta))$ for all $\eta\in S\cap Y_\g^*$.
If $\g$ is a limit ordinal or 0 then put 
$Y_{\g+1}=Y_\g\cup\{\eta\in {^\g\mu}\}$ and there is nothing to check.
If $\g=\d+1$ where $\d$ is a limit ordinal or 0,
let $Y_{\g+1}=Y_\g\cup\bigcup\{Succ_T(\eta):\eta\in Y_\g\cap {^\g\mu}\}$.
Now if $\eta\in S\cap Y_{\g+1}$ for some $\eta\in {^\g\mu}$, then
$\eta|\b\in S\cap Y_\g$ for all $\b<\d$, so $f_0(g(\eta|\b))=f_1(g(\eta|\b))$
for all $\b<\d$.  As both $f_0$ and $f_1$ are order-preserving and
continuous, it follows that $f_0(g(\eta\conc \lan i\ran))=f_1(g(\eta\conc \lan i\ran))$
for $i=0,1$ so our inductive hypothesis is maintained.

Finally, assume $\g=\d+n$, where $\d$ is a limit ordinal or 0 and $n>1$.
Fix $\eta\in Y_\g\cap {^{\d+n-1}\mu}$ and we specify its successors in $Y_{\g+1}$:
\begin{itemize}
\item  If $\eta\not\in S$ or if $Succ_T(\eta)=\emptyset$, then let 
$Succ_{Y_{\g+1}}(\eta)=Succ_T(\eta)$ and there is no problem.

\item  If $\eta\in S$ and $Succ_T(\eta)=\{\eta\conc\lan \al\ran:\al\in C\}$ for some
club subset of an uncountable regular cardinal $\l_\eta$, 
then our hypotheses imply that $f_0(g(\eta))=f_1(g(\eta))$ and
$\{g(\eta\conc\lan\al\ran):\al\in C\}$ is a continuous, strictly increasing (or decreasing)
sequence converging to $g(\eta)$.  Thus, as both $f_0$ and $f_1$ are order-preserving
and continuous, there is a club $C'\sbq C$ such that $f_0(g(\eta\conc\lan\al\ran))=f_1(g(\eta
\conc\lan\al\ran))$ for all $\al\in C'$.  So put $Succ_{Y_{\g+1}}(\eta)=\{\eta\conc\lan\al\ran:\al\in C'\}$.
\endproof

\end{itemize}

As noted above, these invariants are not sufficient to determine the
isomorphism type of an element of $\Khom$.
In particular, the second invariant does not specify which elements
of the representation are in $I$ (as opposed to $\Ibar$).
This affords considerable freedom in choosing the
isomorphism type of the order.  The family of structures in the example below
was first studied by Conway~\cite{Conway} and was later used as an example by Nadel and
Stavi~\cite{NS}.

\begin{example}  \label{many}
%There is a family ${\cal F}$ of $2^{\aleph_1}$
%$L_{\infty,\aleph_1}$-equivalent
%doubly transitive linear orders, each of size
%$\aleph_1$, each sharing the same $\aleph_1$-representation,
%and
%$I\sim I'$ for all
%$I,I'\in {\cal F}$, 
%yet the orders are pairwise non-embeddable.
There is a family  of $2^{\aleph_1}$
locally isomorphic, $L_{\infty,\aleph_1}$-equivalent
doubly transitive linear orders of size
$\aleph_1$, all of whom have 
isomorphic $\aleph_1$-representations;
yet the orders are pairwise non-embeddable.
\end{example}

Let ${\cal S}$ be a collection of $2^{\aleph_1}$ stationary subsets of
$\o_1\setminus\{0\}$ with $X\setminus Y$ stationary for all distinct $X,Y\in {\cal S}$
(see \cite{Solovay} for a construction of such a family).
As notation, let $\QQ^{\ge 0}$ be the set $\QQ\cap [0,\infty)$. For
$X\in {\cal S}$, let 
$$I_X=\sum_{i\in\o_1} J_i^X\qquad \hbox{where}\qquad
 J_i^X=\left\{\begin{array}{ll}
\QQ & \mbox{if $i\not\in X$;}\\
   \QQ^{\ge 0} & \mbox{if $i\in X$.}
\end{array}\right.$$
Clearly $(a,b)\cong \QQ$ for all $a<b$ from $I_X$,
so $I_X\sim I_Y$ for all $X,Y\in {\cal S}$.  
It was first noted by Silver that for any sets $X,Y\in {\cal S}$, the set ${\cal B}(X,Y)$
of all order-preserving partial functions $f:I_X\ra I_Y$, whose domain $D$ is a proper
initial segment of $I_X$ such that $I_X\setminus D$ has no least element, and whose
range $R$ is a proper initial segment of $I_Y$ such that $I_Y\setminus R$ has no least
element, is an $\aleph_1$-back and forth system; hence the orders $I_X$ and $I_Y$
are $L_{\infty,\aleph_1}$-equivalent.
As the Dedekind completions of the $I_X$'s are
isomorphic we can identify them.  After this identification, each of the orders $I_X$
share the same $\aleph_1$-representation, namely
$(T,g)$, where
$T=\{\lan 0\ran,\lan 1\ran\}\cup\{\lan 1,\d\ran:\d\in\o_1\}$
and $g(\lan 1,\d\ran)$ is the element of the Dedekind completion realizing the cut
preceding $J_\d$ for all $\d>0$.

It remains to show that $I_X$ is not embeddable in $I_Y$ whenever $X\neq Y$.
(This was proved in \cite{Conway} but is repeated here for convenience.)
So fix $X\neq Y$ and assume by way of contradiction that there is an
embedding $f:I_X\ra I_Y$.  
It is readily verified that the set
$$C=\{\al\in\omega_1: f(\sum_{i\in\al} J_i^X)=\sum_{i\in\al} J_i^Y\}$$
is a club subset of $\omega_1$.  Thus, since $X\setminus Y$ is stationary, there is
an $\alpha\in C\cap X\setminus Y$.  But $\displaystyle I_X\setminus \sum_{i\in\alpha} J_i^X$
has a least element, whereas 
$\displaystyle I_Y\setminus \sum_{i\in\alpha} J_i^Y$ does not, which is a contradiction.
\endproof

%We will use the existence of a $\mu$-representation to form a desirable partition of
%a doubly transitive linear order.
Despite the limitations demonstrated by  the example above, 
the invariants described in this section do allow us to
obtain an upper bound on the Karp complexity of $\Khom$.
The following definitions
establish our notation.

\begin{defn}  For $D\sbq\Ibar$, a {\it $D$-cut\/} $\nu$ is a partition
of $D$ into two sets, $D^-_\nu$ and $D^+_\nu$ (either may be empty)
such that $D^-_\nu\cup D^+_\nu=D$, $D^-_\nu\cap D^+_\nu=\e$, and
$D^-_\nu$ is downward closed.  We write $\nu=(D^-_\nu,D^+_\nu)$
and let $I(\nu)=\{x\in I:D^-_\nu<x<D^+_\nu\}$. 
%$\d^-(\nu)=\sup(D^-_\nu)$.
\end{defn}

\begin{defn}
Suppose $I$ and $J$ are two linear orders.  
If $D\sbq\Ibar$ and $f:D\ra\Jbar$ is any order-preserving function then
$f(\nu)$ is the $f(D)$-cut $(f(D^-_\nu),f(D^+_\nu))$.
A function $f:D\ra\Jbar$ is {\it proper\/} if $\{-\infty,+\infty\}\sbq D$ and $f$
is order-preserving,
continuous, $f(-\infty)=-\infty$, $f(+\infty)=+\infty$,
and satisfies $d\in I\Leftrightarrow f(d)\in J$ for all $d\in D$.
\end{defn}

If $D\sbq \Ibar$ and $f:D\ra \Jbar$ is a proper function, then $I\setminus D$ and $J\setminus f(D)$
are partitioned into corresponding families of $D$-cuts and $f(D)$-cuts.  
The following definitions measure the  similarity of these cuts.
%We extend our
%defintion of Rank from Definition~\ref{Rank} to the class  of proper functions as follows.

\begin{defn}  Two (possibly empty)
linear orders $I$ and $J$ are {\it $(\mu^+,\alpha)$-equivalent,\/}
written $I\equiv_{\mu^+,\al} J$, if $I$ and $J$ are elementarily equivalent and
the empty function in $\FF_{\mu^+}(I,J)$ has
Rank at least $\alpha$ (see Definition~\ref{Rank}).
\end{defn}

By allowing linear orders to be empty and by insisting on elementary equivalence
we intend that
$I=\emptyset$ if and only if $J=\emptyset$ and $|I|=1$ if and only if $|J|=1$
whenever $I\equiv_{\mu+,\al} J$ for some ordinal $\alpha$.

\begin{defn}  If $D\sbq \Ibar$ and $f:D\ra\Jbar$ is proper, then $f$ is {\it $\alpha$-strong\/}
if $I(\nu)\equiv_{\mu^+,\al} J(f(\nu))$ for all $D$-cuts $\nu$.
\end{defn}

%It is easily verified by induction on $\al$ that for $\al\ge 1$,
%a function $f\in \FF_{\mu+}(I,J)$ has Rank at least $\al$ if and only if 
%the canonical extension $\tilde{f}$ of $f$ is $\al$-strong.
%
%
%
%It is easy to see that if $f:A\ra I_1$ has Rank at least 2, then many pathological cases
%are eliminated.  In particular, any such $f$ is
%continuous and has a unique 
%extension to a proper map $\tilde{f}:\lim(A)\cup\{-\infty,+\infty\}\ra\Ibar_1$.
%Further, $I(\nu)$ is elementarily equivalent to $I_1(f(\nu))$ for any 
%such function $f$.
%The following lemma is easily established by induction on $\alpha$.

If $f\in\FF_{\mu^+}(I,J)$ has domain $A$ and has Rank at least 2,
then it is easily seen that $f$ is continuous and extends uniquely to a proper function 
$$g:A\cup\lim(A)\cup\{-\infty,+\infty\}\ra\Jbar,$$
where $\lim(A)$ denotes the set of limit points of $A$ in $\Ibar$.
Also, it is easily established by induction on $\al\ge 1$ that if 
$g:D\ra\Jbar$ is a proper function with domain $D\sbq\Ibar$ and the restriction $f=g|(D\cap I)$ is 
in $\FF_{\mu^+}(I,J)$, then $g$ is $\al$-strong if and only if Rank$(f)\ge\al$.

For $\al\ge 1$ the class of $\al$-strong proper functions has desirable closure properties.
It is routine to show that the restriction of any $\al$-strong proper function to
any set that contains $\{-\infty,+\infty\}$ is also proper and $\al$-strong.
As well, we have the following lemma, which is proved by a straightforward induction on $\alpha$.

\begin{lemma}  \label{partition}
Let $\al\ge 1$.  Suppose that 
$D\sbq \Ibar$, $f:D\ra\Jbar$ is an $\alpha$-strong proper function,
and for each $D$-cut $\nu$ there is a set $E_\nu\sbq \overline{I(\nu)}$ and an
$\al$-strong proper function $g_\nu:E_\nu\ra\overline{J(f(\nu))}$.
Then $f\cup\bigcup_\nu g_\nu$ is proper and $\al$-strong.
\end{lemma}

%\begin{lemma}  \label{partition} Fix $\al\ge 2$.
%%If $f:D\ra\Ibar_1$ is proper and for every $D$-cut
%$\nu$ with $I_0(\nu)$ nonempty, there is
%a subset $A_\nu\sbq I_0(\nu)$ and an $\al$-strong
%map $g_\nu:A_\nu\ra I_1(f(\nu))$, then the
%combined map $h:(D\cap I_0)\cup\bigcup_\nu A_\nu\ra I_1$
%defined by
%$$h(a)=\left\{\begin{array}{ll}f(a) & \mbox{if $a\in D\cap I_0$;}\\
%g_\nu(a) & \mbox{if $a\in I_0(\nu)$}\end{array}
%\right.$$
%is $\al$-strong.
%\end{lemma}
%

\begin{lemma}  \label{cases}
Let $I_0,I_1\in\Khom$ satisfy $I_0\sim I_1$ and $I_0\equiv_{\mu^+,\al} I_1$ for some ordinal 
$\al\ge 2$.
Assume that $A\sbq I_0$ is of size at most $\mu$ and satisfies
%$$\big[A\ \hbox{is bounded below or}\ \coi(A)=\ao\big]\ 
%\hbox{and}\ \big[A\ \hbox{is bounded above or}\
%\cof(A)=\ao\big]$$
\begin{enumerate}
\item $A$ is bounded below or $\coi(A)=\ao$; and
\item $A$ is bounded above or $\cof(A)=\ao$.
\end{enumerate}
Then there is an $f\in\FF_{\mu^+}(I_0,I_1)$ with domain $A$ of Rank at least $\alpha$.
\end{lemma}

\bp  
%First, suppose $A$ is bounded both above and below.
%Choose $a,b\in I_0$ with $a<A<b$, choose any $c<d$ from 
%$I_1$, and choose an isomorphism $f:[a,b]\ra[c,d]$.
%Then clearly $f\r A$ is $\al$-strong.
%Next, suppose  $A$ is unbounded both above and below
%and $\coi(A)=\cof(A)=\ao$.
%We claim that $I_0\cong I_1$ in this case.
%Since $I_0\equiv_{\mu^+,\al} I_1$, it follows that
%$\coi(I_1)=\cof(I_1)=\ao$.  Choose strictly increasing sequences
%$\lan a_n:n\in {\bf Z}\ran$ from $A$ and
%$\lan c_n:n\in {\bf Z}\ran$ from $I_1$,
%each 
%both cofinal and coinitial in $I_0$ and $I_1$, respectively.
%For each $n\in {\bf Z}$ there is an isomorphism $f_n:[a_n,a_{n+1}]
%\ra[c_n,c_{n+1}]$, so their union is an isomorphism between
%$I_0$ and $I_1$.
%The hybrid cases (e.g., $A$ bounded below and $\cof(A)=\ao$)
%are handled similarly.
We show that in fact $A$ is contained in an interval of $I_0$
which is isomorphic to an interval of $I_1$.  This interval will be of the form
$(a,b)$, where $a$ is a lower bound for $A$ if one exists, or the symbol $-\infty$,
and $b$ is defined similarly.  Take as a typical case that in which $a\in I_0$ and $b=\infty$.
Then we claim that the interval $(a,\infty)$ is isomorphic to $(a',\infty)$ for any
$a'\in I_1$.  The point is that $(a,\infty)$ has cofinality $\aleph_0$,
hence $(a',\infty)$ does by $(\mu^+,\al)$-equivalence.
So we can build the desired isomorphism in a countable sequence of steps, using
double transitivity and the local isomorphism of $I_0$ and $I_1$.

As well, it follows from the relations
$I_0\sim I_1$ and $I_0\equiv_{\mu^+,\al} I_1$
and another instance of double transitivity that 
the intervals $(-\infty,a)$ and $(-\infty,a')$ are $(\mu^+,\al)$-equivalent.
Thus, the 
the restriction of the isomorphism to $A$ has Rank at least $\al$.
\endproof

The following Proposition is the key to the proof of Theorem~\ref{upper}.
Before embarking on it, we introduce some more notation.
For $C\sbq \l$, let 
$$C_0=\{\al\in C:\al\ \hbox{is a limit point of}\ C\cap\al\}$$
and for $T$ a $\mu$-decomposition tree, let $T_0$ be the club subtree of
$T$ satisfying $Succ_{T_0}(\eta)=\{\eta\conc\lan\al\ran:\al\in C_0\}$,
where $Succ_T(\eta)=\{\eta\conc\lan\al\ran:\al\in C\}$ for all non-maximal nodes $\eta\in T^*_0$.
Note that if $(T,g)$ is a $\mu$-representation of $\Ibar$, then
$(T_0,g|T^*_0)$ is also a $\mu$-representation of $\Ibar$ with the additional property
that
$g(\eta)$ either has cofinality or coinitiality at most $\mu$
for all $\eta\in T_0^*\setminus\{\langle 0\rangle,\langle 1\rangle\}$.

\begin{prop}  \label{key}
Assume $I_0,I_1\in\Khom$, $I_0\sim I_1$ and $I_0\equiv_{\mu^+,\o} I_1$.
If $A\sbq I_0$ and $|A|\le\mu$, then there is a function
$f:A\ra I_1$ of Rank at least $\omega$.
\end{prop}

\bp
Pick $A\sbq I_0$ of size at most $\mu$.
In order to produce a $h:A\ra I_1$
of Rank at least $\o$,
we first construct a desirable proper function $j:D\ra\Ibar_1$.
Choose a $\mu$-representation $(T,g)$ of $\Ibar_0$.  By passing to the
subtree $T_0$ in the notation preceding this proposition, we may assume that
$g(\eta)$ either has cofinality or coinitiality at most $\mu$ for all 
$\eta\in T^*\setminus\{\langle 0\rangle,\langle 1\rangle\}$.
Let $B=B_L\cup B_R$, where
\begin{eqnarray}
B_L & =  &  \{\eta\in T^*:dir(\eta)=\hbox{LEFT and $A$ is cofinal in $I_\eta$}\}\quad
\hbox{and}  \nonumber \\
B_R & =  &  \{\eta\in T^*:dir(\eta)=\hbox{RIGHT and $A$ is coinitial in $I_\eta$}\}.  \nonumber
\end{eqnarray}
We claim that $B$ has size at most $\mu$.    To see this, it suffices by symmetry to
show that $|B_L|\le\mu$.  Recall that for every successor ordinal $\al$,
the intervals 
$\{I_\eta:\eta\in T^*\cap {^{\al}\mu}\}$ are disjoint.
Since $\eta\in B_L$ implies $A\cap I_\eta\neq\emptyset$,
this implies $|B_L\cap {^\al\mu}|\le\mu$
for all successor ordinals $\alpha$.
Further, since $|A|\le\mu$, we can choose a successor ordinal
$\g$ so that for every pair $a,a'\in A$, there is $\eta\in T^*$ of length
less than $\g$ satisfying $a<g(\eta)<a'$ whenever there is any
$\nu\in T^*$ with $a<g(\nu)<a'$.
But now, by our choice of $\g$, if $\nu,\nu'\in B_L$
have length $>\g$ and have $\nu|\g=\nu'|\g$,
then $$g(\nu)=\sup(A\cap I_{\nu|\g})=g(\nu'),$$
so $\nu=\nu'$ and $I_{\nu|\g}\cap A\neq\emptyset$.
Thus,
$$|\{\nu\in B_L:lg(\nu)>\g\}|\le|\{\eta\in T^*:lg(\eta)=\g
\ \hbox{and}\ I_\eta\cap A\neq\emptyset\}|\le\mu$$
so $|B_L|\le\mu$.

Let $B'\contains B$ be a closed subset of $T$  of size at most $\mu$.
As $g(\eta)$ has cofinality or coinitiality at most $\mu$ in $\Ibar_0$
for each $\eta\in B'\setminus\{\langle 0\rangle,\langle 1\rangle\}$, 
there is a set $X\sbq I_0$ of size at most $\mu$
such that $g(B')\sbq\lim(X)\cup \{-\infty,+\infty\}$.
Since $I_0\equiv_{\mu^+,\o} I_1$, for each $n\ge 2$ we can choose an
order-preserving
$j_n:X\ra I_1$ of Rank at least $n$.  As $g(B')\sbq\lim(X)$, each $j_n$ extends uniquely
to a proper function (also called $j_n$) from $X\cup g(B')$
to $\Ibar_1$.  
As $B'\sbq T$ is closed, by Lemma~\ref{closed}
there is a club subtree $T'_n$ for each $n\ge 2$ such that 
$j_n(g(\eta))=j_{n+1}(g(\eta))$ for all $\eta\in B'\cap T'_n$.
Let $Y=\bigcap_{n\ge 2} T'_n$
and let $D=\{g(\eta):\eta\in B'\cap Y_0\}$, where $Y_0$ is the
club subtree of $Y$ described in the notation preceding this proposition.
As the functions $j_n$ agree on $D$ for all $n\ge 2$, we let
$j:D\ra\Ibar_1$ denote this common (proper) function.
As each $j_n$ was $n$-strong, the function $j$ is $\o$-strong.

By Lemma~\ref{partition}, in order to
ascertain the existence of an $\o$-strong $h:A\ra I_1$, it suffices
to construct an order-preserving function $f:A\cap I_0(\nu)\ra I_1(j(\nu))$ 
of Rank at least $\o$ for
every $D$-cut $\nu$ of $I_0$.
So fix a $D$-cut $\nu=(D^-_\nu,D^+_\nu)$.  
We finish the proof by showing that the hypotheses of Lemma~\ref{cases}
are satisfied for $I_0(\nu)$ and $I_1(j(\nu))$.
As $I_0(\nu)$ and $I_1(j(\nu))$ are convex subsets of $I_0$ and $I_1$
respectively, $I_0(\nu)\sim I_1(j(\nu))$.
Since $j$ is $\o$-strong,  $I_0(\nu)\equiv_{\mu^+,\o} I_1(j(\nu))$.
Finally, assume by way of contradiction that $A\cap I_0(\nu)$ is unbounded above in $I_0(\nu)$
and has uncountable cofinality.  
(The case of $A\cap I_0(\nu)$ unbounded below in $I_0(\nu)$ of uncountable cardinality 
is symmetric.)
Let $b=\sup(A\cap I_0(\nu))$ and let $\k=\cof(A\cap I_0(\nu))$.
We will obtain a contradiction by showing that $b=\sup(D^-_\nu)$,
which would make $I_0(\nu)$ empty.
First, since $Y_0$ is a club subtree of $T$ and $b=\inf(D^+_\nu)$,
$b=g(\eta)$ for some $\eta\in Y_0$.  As we assumed $A$ cofinal below $b$,
$b\in B$ as well.
There are now four cases to consider, all of which imply $b=\sup(D^-_\nu)$ or
contradict our hypotheses.
\smallskip\par\noindent
{\bf Case 1.}  $dir(\eta)=$RIGHT and $lg(\eta)=\d+1$ where $\d$ is a limit ordinal or 0.
\smallskip
\par
Say $\eta=\rho\conc \lan 0\ran$.  Since $\cof(b)=\k>\ao$
there is a strictly increasing sequence of limit ordinals $\lan\g_i:i<\k\ran$
such that $b=\sup\{g(\rho|(\g_i+1)):i<\k\}$.  Since $B'$ is closed,
$\rho|\g\in B'$ for all $\g<lg(\rho)$, so $g(\rho|\gamma)\in D_\nu$ and
$b=\sup(D^-_\nu)$.
\smallskip\par\noindent
{\bf Case 2.}  $dir(\eta)=$RIGHT and $lg(\eta)=\d+n$ for some $n>1$.
\smallskip\par
Say $\eta=\rho\conc\lan\al\ran$ for some $\al\in C_0$, where 
$C$ is such that $Succ_Y(\rho)=\{\rho\conc\lan\b\ran:\b\in C\}$.
As $\cof(b)=\k$ there is a continuous, strictly increasing sequence of ordinals
$\lan\b_i:i<\k\ran$ from $C$ with limit $\al$.
Again, as $B'$ is closed, 
$\rho\conc\lan \b_i\ran\in B'$
for all $i\in \k$.  It follows that $\rho\conc\lan\b_i\ran\in B'\cap Y_0$
for all limit ordinals $i\in \k$, so again $b=\sup(D^-_\nu)$.
\smallskip\par\noindent
{\bf Case 3.}
$dir(\eta)=$LEFT and $\eta$ is not maximal in $Y_0$.
\smallskip\par
Say $Succ_{Y_0}(\eta)=\{\eta\conc\lan\al\ran:\al\in C_0\}$.
As $A$ is unbounded below $b$ and $\k>\ao$, there is a club $C'\sbq C_0$
such that $A$ is unbounded below
$g(\eta\conc\lan\al\ran)$ for all $\al\in C'$.
Thus, $\eta\conc\lan\al\ran\in B\cap Y^*_0$ for all $\al\in C'$, so 
again $b=\sup(D^-_\nu)$.
\smallskip\par\noindent
{\bf Case 4.}
$dir(\eta)=$LEFT and $\eta$ is maximal in $Y_0$.
\smallskip
\par
As $\eta$ maximal in $Y_0$ implies $\eta$ maximal in $T$,
it follows from the definition of a $\mu$-representation 
that $\cof(g(\eta))=\ao$ or $\cof(g(\eta))>\mu$.
However, we assumed that $\cof(g(\eta))>\ao$ and 
$A$ witnesses that $\cof(g(\eta))\le\mu$, so both are impossible.
\endproof

Our theorem now follows easily.

\begin{theorem}  \label{upper}
$KC_{\mu^+}(\Khom)\le\o$ for all  uncountable cardinals $\mu$.
\end{theorem}

\bp
Fix $I\in\Khom$ and an uncountable cardinal $\mu$.  Let $f\in \FF_{\mu^+}(I)$ have Rank
at least $\o$.  We claim that Rank$(f)\ge\o+1$.
To see this, it suffices by
symmetry to show that if
$A\sbq I$, $|A|\le\mu$  
then there is a function $g\in \FF_{\mu^+}(I)$ extending $f$ of Rank at least $\o$
with $A\sbq \dom(g)$.  
So fix such a set $A$ and let $\tilde{f}$ denote the proper function extending $f$
with domain $\dom(f)\cup\{-\infty,+\infty\}$.
Since Rank$(f)\ge\o$, $\tilde{f}$ is $\o$-strong.
Now fix a 
$\dom(\tilde{f})$-cut $\nu$.
Clearly, 
$I(\nu)\equiv_{\mu^+,\o} I(f(\nu))$ and $I(\nu)\sim I(f(\nu))$, so
it follows from Proposition~\ref{key} that there is a  function $g_\nu:
A\cap I(\nu)\ra I(\tilde{f}(\nu))$ in $\FF_{\mu^+}(I(\nu),I(\tilde{f}(\nu))$
of Rank at least $\o$.
Thus, it follows from Lemma~\ref{partition} that  the proper function 
$g=\tilde{f}\cup\bigcup\{g_\nu:\nu$ a $\dom(\tilde{f})$-cut$\}$
is $\o$-strong, hence the restriction of $g$ to $A\cup\dom(f)$ has Rank at least $\o$.
\endproof

\def\nubar{\bar{\nu}}
\def\betabar{\bar{\b}}
\section{The $\o$-independence property}  \label{omegaindep}

This section is devoted to proving that any pseudo-elementary class with the 
$\o$-independence property (see Definition~\ref{o-indep}) is not controlled.
We begin the section by proving Proposition~\ref{criterion},
which will provide us with
a method for concluding that $KC_\l(\KK)=\infty$
by looking at the family of  $\l$-partial isomorphisms from one element of $\KK$ into
another.

\begin{defn}  An {\it $\o$-tree\/}
$\T$ is a downward closed subset of $^{<\o}\l$ for some ordinal $\l$.
We call $\T$ {\it well-founded\/} if it does not have an infinite branch.
For a tree $\T$ and $\eta\in \T$, the {\it depth of $\T$ above $\eta$\/}, $dp_{\T}(\eta)$
is defined inductively by
$$dp_{\T}(\eta)=\left\{\begin{array}{ll}
\sup\{dp_{\T}(\nu)+1\}:\eta\lhd\nu\}&\mbox{if $\eta$ has a successor}\\
                           0&\mbox{otherwise.}
\end{array}\right. $$

and the depth of $\T$, $dp(\T)=dp_{\T}(\langle\rangle)$.
\end{defn}

Clearly, $dp(\T)<\infty$ if and only if $\T$ is well-founded.
The most insightful example is that for any ordinal $\al$, the tree $des(\al)$ consisting of
all descending sequences of ordinals $<\alpha$ ordered by initial segment
has depth $\alpha$.
The proof of the following lemma is reminiscent of the proof of Morley's Omitting
Types Theorem.

\begin{lemma}  \label{tree} 
If $\T\sbq {^{<\o}\l}$ is well-founded and $dp(\T)\ge\k^+$, then for 
any coloring $c:\T\ra \k$, there is a subtree $\S\sbq \T$ of depth at least $\o$
such that $c|_{\S\cap {^n\l}}$ is constant for each $n\in\o$.
\end{lemma}

\bp
Since $dp_{\T}\ge\k^+$, there is an $\eta\in\T$ with $dp_\T(\eta)=\k^+$.
Thus, by concentrating on subtrees extending $\eta$, we may assume that
$dp(\T)=\k^+$.

For each $n\in\o$ we will name a subset $X_n\sbq \k^+$ of size $\k^+$
and a function $f_n:X_n\ra \T\cap {^n\l}$ such that $X_{n+1}\sbq X_n$,
every element of $f_{n+1}(X_{n+1})$ is a successor of an element of $f_n(X_n)$,
$dp_\T(f_n(\al))\ge\al$ and $c|_{f_n(X_n)}$ is constant.

To begin, let $X_0=\k^+$ and let $f_0:X_0\ra \{\langle\rangle\}$.
Given $X_n$ and $f_n$ satisfying our demands, we
define  $X_{n+1}$ and $f_{n+1}:X_{n+1}\ra \T\cap {^{n+1}\l}$ as follows.
For $\al\in X_n$, let $\b$ be the least element of $X_n$ greater than $\al$.
As $dp_{\T}(f_n(\b))\ge\b$, we can define $f_{n+1}(\al)$ to be a successor of
$f_n(\b)$ of depth at least $\al$.  Since $X_n$ has size $\k^+$, let $X_{n+1}$
be a subset of $X_n$ of size $\k^+$ such that $c|_{f_{n+1}(X_{n+1})}$ is monochromatic.

Now let $R=\{f_n(\b_n):n\in\o\}$, where $\b_n$ is the least element of $X_n$
and let $\S$ be the subtree of $\T$ generated by $R$.
\endproof

Suppose that $N\equiv M$ and $\emptyset=A_0\sbq A_1\sbq\dots\sbq N$ is an $\o$-sequence
of subsets of $N$ of size less than $\l$.  Let
$$\T_n=\{{\rm Range}(f):f\in\FF_\l(N,M), \ f\ \hbox{has domain $A_n$}\}$$
and let $\T=\bigcup\{\T_n:n\in\o\}$ be a tree under inclusion.
Typically $\T$ will be an $\o$-tree and
we can ask whether or not it is well-founded.  The relationship
between this question and Karp complexity is partially explained by the following
proposition.

\begin{prop}  \label{criterion}
If $KC_\l(\KK)<\infty$ then there is an ordinal
$\al^*$ such that whenever $N\equiv M\in\KK$ and 
$\emptyset=A_0\sbq A_1\sbq\dots\sbq N$ are chosen with $|A_i|<\l$,
then the induced tree $\T$ either has depth at most $\al^*$ or
has an infinite branch.
\end{prop}

\bp  If $KC_\l(\KK)<\infty$ then by Proposition~\ref{equiv}, there is a cardinal
$\k$ bounding the number of $L_{\infty,\l}$-types realized in elements of $\KK$.
We claim that $\al^*=\k^+$ has the desired property.  To see this, choose
$N\equiv M$ from $\KK$ and 
$\emptyset=A_0\sbq A_1\sbq\dots\sbq N$ and assume that $dp(\T)\ge\k^+$.
By Lemma~\ref{tree}, there is a subtree $\S$ of $\T$ of depth $\o$
such that the $L_{\infty,\l}$-types of the elements of $\S$ depend only
on their level in $\S$.
In particular, for each $n$ there is an element $B_n\in\S$ at level $n$
that has a successor in $\S$.  Consequently, for each $n\in\o$
the $L_{\infty,\l}$-formula
$$\Theta(X_n)=\exists Y_n tp_{\infty,\l}(X_n,Y_n)=tp_{\infty,\l}(B_{n+1})$$
is implied by $tp_{\infty,\l}(B_n)$.
Applying this iteratively produces an elementary partial function $f:N\ra M$ with
domain $\bigcup\{A_n:n\in\o\}$, so $\T$ has an infinite branch.\endproof

\begin{defn}  \label{o-indep}
A class $\KK$ of $L$-structures has the {\it $\o$-independence property\/}
if there is a set $\{\phi_n(\xbar_0,\dots,\xbar_{n-1},\ybar_n):n\in\o\}$
of $L$-formulas such that for all $M\in\KK$ there
is a sequence $\lan \abar_i:i<\o\ran$ from $M$ such that for all $n\in\o$
and all functions 
$f:n\ra \{0,1\}$
there is a sequence $\lan \bbar_i:i<n\ran$ from $M$ such that for all $i<n$,
$$M\models\phi_i(\bbar_0,\dots,\bbar_{i-1},\abar_i) \qquad \hbox{if and only if}
\qquad f(i)=1.$$
\end{defn}

As an example, the model completion of the empty theory in the language
$L=\{R_n:n\in\o\}$ consisting of one $n$-ary relation for every $n$
is a complete, simple theory with the $\o$-independence property.
(In this example, the $\ybar_n$'s do not appear.)  Clearly, if
$\KK$ has the $\o$-independence property, then $\KK$ has the independence
property.  However, the theory of the random graph has the independence
property, but fails to have the $\o$-independence property.
We remark that despite this failure,
the theory of the random graph is not controlled. We do
not attempt to prove this assertion here.

Our interest in the notion of $\o$-independence 
is largely captured by the proposition given below.

\begin{defn}  An {\it ordered multigraph\/} is a structure $(G,<,R_n)_{n\in \o}$
where $<$ is interpreted as a linear order and each $R_n$ is a symmetric 
$n$-ary relation on $G$.
\end{defn}

\begin{prop}  \label{coding}
If $L_1\contains L_0$, $T_1$ is an $L_1$-theory
with Skolem functions and $\KK$, the
class of reducts of models of $T_1$ to $L_0$ has the $\o$-independence
property witnessed by $\{\phi_n:n\in\o\}$
then for every ordered multigraph $(G,<,R_n)_{n\in\o}$
there is a structure $M_G\in\KK$ and sequences $\lan \abar_n:n\in\o\ran$
and $\lan \bbar_g:g\in G\ran$ from $M_G$
such that
\begin{enumerate}
\item  $M_G$ is the $L_1$-Skolem hull of $\{\abar_n:n\in\o\}\cup\{\bbar_g:g\in G\}$;
\item If $g_1,\dots,g_n$ and $h_1,\dots,h_n$ have the same quantifier-free type in
$(G,<,R_n)_{n\in\o}$ then the sequences 
$\bbar_{g_1},\dots,\bbar_{g_n}$ and 
$\bbar_{h_1},\dots,\bbar_{h_n}$ have the same type over $\{\abar_n:n\in\o\}$ in $M_G$;
\item  $M_G\models \phi_n(\bbar_{g_1},\dots,\bbar_{g_n},\abar_n)$ if and only if
$G\models R_n(g_1,\dots,g_n)$ for all $n$ and all $g_1,\dots,g_n$ from $G$.
\end{enumerate}
\end{prop}

The proof of Proposition~\ref{coding} is word for word like the proof of the existence
of Ehrenfeucht-Mostowski models for unstable pseudo-elementary classes 
(see e.g., Section~11.3 of \cite{Hodges})
but with the Ne\u{s}et\u{r}il-R\"odl theorem (see \cite{Nesetril}
or \cite{AbLeo})
in place of Ramsey's theorem.

The following lemma tells us that we need not explicitly
consider the constants $\{\abar_n:n\in\o\}$ in the proof of Theorem~\ref{uncontrolled}.

\begin{lemma}  \label{constants}
Let $\KK$ be a class of $L$-structures
and let $C$ be a set of fewer than $\l$ new constant symbols.
Let $\KK^*$ be the class of all expansions of
elements of $\KK$ to $L\cup C$-structures.
Then $KC_\l(\KK^*)\le KC_\l(\KK)$.
\end{lemma}

\bp  For any $M^*\in\KK^*$, let $M$ be its reduct to the language of $L$.
For every partial function $f\in\FF_\l(M^*)$, let $\ff\in\FF_\l(M)$
be the extension of $f$ that is the identity on every element of $C^M$.
It is easy to show by induction that Rank$_{\FF^*}(f)={\rm Rank}_\FF(\ff)$.
Hence, $KC_\l(M^*)\le KC_\l(M)$, so $KC_\l(\KK^*)\le KC_\l(\KK)$.\endproof

The other theorem we will need is that there exist very complicated colorings of a number of
cardinals.  As notation, for $x$ a finite subset of $\mu$,
let $x^m$ denote the $m^{{\rm th}}$ element of $x$ in 
increasing order.  Following the notation in \cite{Sh:g},
let $Pr_0(\mu,\mu,\ao,\ao)$
denote the following statement:
\begin{itemize}
\item  There is a symmetric two-place function $c:\mu\times\mu\ra \o$ such that
for every $n\in\o$, every collection of $\mu$ disjoint, $n$-element subsets 
$\{x_\al:\al\in\mu\}$ of $\mu$,
and every function $f:n\times n\ra \o$,
there are $\al<\b<\mu$ such that 
$$c(x^m_\al,x^{m'}_\b)=f(m,m')$$
for all $m,m'<n$.
\end{itemize}

It is shown in \cite{Sh:327} that $Pr_0(\l,\l,\ao,\ao)$ holds
for an uncountable cardinal $\l$ whenever there exists a nonreflecting
stationary subset of $\l$.  (A stationary subset $S\sbq \l$ is nonreflecting
if $S\cap\al$ is not stationary in $\al$ for all limit ordinals $\al<\l$.)
In particular, $Pr_0(\aleph_3,\aleph_3,\ao,\ao)$ holds.
More recently, in \cite{Sh:572}
the second author has shown that $Pr_0(\aleph_2,\aleph_2,\ao,\ao)$
holds as well.  This suffices for our purpose.
See \cite{Sh:g} for more of the history of $Pr_0$ and its cousins.

The following Lemma recasts
$Pr_0(\mu,\mu,\ao,\ao)$ into the form we will use in 
the proof of Theorem~\ref{uncontrolled}.

\begin{lemma}  \label{useful}
Let $c:[\mu]^2\ra\o$ witness $Pr_0(\mu,\mu,\ao,\ao)$.
For every $k,n\in\o$, every collection $\{x_\al:\al\in\mu\}$
of $\mu$ disjoint, $n$-element subsets of $\mu$, and every family 
of colorings $\{f_{i,j}:n^2\ra\o:i<j<k\}$, there are
$\b_0<\b_1<\dots<\b_{k-1}$ such that
$$c(x^m_{\b_i},x^{m'}_{\b_j})=f_{i,j}(m,m')$$
for all $i<j<k$ and all $m,m'<n$.  
\end{lemma}

\bp  
Fix $k,n$, $\{x_\al:\al\in\mu\}$,
and $\{f_{i,j}:i<j<k\}$ satisfying the hypotheses.
Without loss, we may assume that $x_\al^{n-1}<x^0_{\al+1}$ for all $\al$.
For $\al$ limit, let $y_\al=\bigcup\{x_{\al+i}:i<k\}$
and let $W_0=\{\al\in\mu:\al\ {\rm limit}\}$.
By induction on $k'\le k$ we will build a sequence
$\b_0<\b_1<\dots<\b_{k'-1}$ and a subset $W_{k'}$
of size $\mu$ such that 
$W_{k'+1}\sbq W_{k'}$ and $c(y^{ni+m}_{\b_i},y^{nj+m'}_{\g})=f_{i,j}(m,m')$
for all $m,m'<n$, all $i<j<k$ with $i<k'$, and all $\g\in W_{k'}$.
For $k'=0$ there is nothing to do.  Assuming $\b_0<\dots<\b_{k'-1}$
and $W_{k'}$ have been chosen,
it follows from $Pr_0(\mu,\mu,\ao,\ao)$
that there is $\b_{k'}$ such that the set
$$\{\g\in W_{k'}:\g>\b_{k'}\
{\rm and}\ c(y^{nk'+m}_{\b_{k'}},y^{jn+m'}_\g)=f_{k',j}(m,m')\ \hbox{for $j>k'$}\}$$
has size $\mu$, hence is a suitable choice for $W_{k'+1}$.
(If there were no such $\b_{k'}$ then one could successively
build a subset $Z$ of $W_{k'}$ of size $\mu$
on which there would be no $\al<\b$ from $Z$ 
satisfying the coloring.)\endproof 

\begin{theorem}  \label{uncontrolled}
Let
$L_1\contains L_0$ be first order languages,
let $T_1$ be an $L_1$-theory and let $\KK$
denote the class of reducts of models of $T_1$ to $L_0$.
If $\KK$ has the $\o$-independence property then $\KK$ is not
controlled.  More precisely, if a cardinal $\mu>|T_1|$ 
is regular and there is a coloring of $[\mu]^2$ satisfying
$Pr_0(\mu,\mu,\ao,\ao)$,
then $KC_{\l}(\KK)=\infty$ for all cardinals $\l>\mu$.
\end{theorem}

\bp  First, by adding countably many constants to the language
$L_0$ and invoking Lemma~\ref{constants},
we may assume that the $\o$-independence of
$\KK$ is witnessed by formulas $\phi_n(\xbar_0,\dots,\xbar_{n-1})$
with no additional constants.
Second, by considering $M^{\rm eq}$ in place of $M$ for each $M\in\KK$,
we may assume that each $\xbar$ is a singleton.
Third, by expanding $T_1$ if necessary, we may assume that it has built-in
Skolem functions.
Fix a coloring
$c:[\mu]^2\ra\o$ that witnesses $Pr_0(\mu,\mu,\ao,\ao)$ and fix an ordinal $\al^*$.
We will 
use the coloring 
to define two rather complicated ordered multigraphs $I$ and $J$
and then
use Proposition~\ref{coding} to get Ehrenfeucht-Mostowski models $M,N\in\KK$
that are built from $I$ and $J$ respectively.  We will find a tree
of $\l$-partial isomorphisms from $N$ into $M$ that is well-founded, yet has
depth at least $\al^*$. Since $\al^*$ was arbitrary,
it follows immediately from 
Proposition~\ref{criterion} that 
$KC_\l(\KK)=\infty$.
So, let
$$des(\al^*)=\{\hbox{strictly decreasing sequences of ordinals}\ <\al^*\}$$
and let $(I,<)$ be the linear order with universe $\mu\times des(\al^*)$,
ordered lexicographically.  Let $(J,<)$ be the linear order
with universe $\mu\times\{\rho_n:n\in\o\}$, where $\rho_n=\lan 0,-1,-2,\dots,-n+1\ran$,
also ordered lexicographically.

As notation, for finite sequences $\eta,\nu$
we write $\eta\lhd\nu$ when $\eta$ is a proper initial segment of $\nu$.
For $t\in I\cup J$, let $t=(\al^t,\eta^t)$, where $\al^t\in\mu$
and $\eta^t$ is a finite, decreasing sequence.
For $s,t\in I\cup J$, we write $s\lhd^* t$ when $\eta^s\lhd\eta^t$.
Fix, for the whole of this section, a partition of $\o\setminus\{0\}$
into disjoint, infinite sets $\{Z_n:n\in\o\}$.

We expand $(I,<)$ into an ordered multigraph 
$(I,<,R_n)_{n\in\o}$
as follows:
We posit that $R_0$ holds, $R_1(t)$ holds for all $t\in I$, and
for $n>1$,
$R_n(t_0,\dots,t_{n-1})$ holds  
if and only if for some permutation $\sigma\in Sym(n)$,
\begin{itemize}
\item $\eta^{t_{\sigma(0)}}\lhd \dots\lhd\eta^{t_{\sigma(n-1)}}$;
\item $lg(\eta^{t_{\sigma(i)}})=i$ for all $i$; 
\item $\al^{t_i}\neq\al^{t_j}$ and $c(\al^{t_i},\al^{t_{j}})\in Z_n$ for all $i<j<n$; and 
\item $c(\al^{t_i},\al^{t_j})=c(\al^{t_k},\al^{t_l})$ for all $i,j,k,l<n$ with $i\neq j$ and $k\neq l$.
\end{itemize}

Similarly, expand $(J,<)$ to an ordered multigraph 
$(J,<,R_n)_{n\in\o}$
by positing that  
$R_0$ holds, $R_1(t)$ holds for all $t\in J$, and
for all $n>1$
$R_n(t_0,\dots,t_{n-1})$ holds 
if and only if for some $\sigma\in Sym(n)$,
\begin{itemize}
\item $\eta^{t_{\sigma(i)}}=\rho_i$ for all $i<n$;
\item $\al^{t_i}\neq \al^{t_j}$ and $c(\al^{t_i},\al^{t_{j}})\in Z_n$ for all $i<j<n$; and 
\item $c(\al^{t_i},\al^{t_j})=c(\al^{t_k},\al^{t_l})$ for all $i,j,k,l<n$ with $i\neq j$ and $k\neq l$.
\end{itemize}

Now build Ehrenfeucht-Mostowski models
$M,N\in\KK$ from $I$ and $J$ respectively that satisfy
Conditions 1-3 of  Proposition~\ref{coding}.  To avoid wanton use of 
nested subscripts, we identify the elements $b_g\in M$ and $g\in I$
(and similarly for $N$ and $J$).

For each $n\in\o$ let $A_n=\{t\in N:lg(\eta^t)<n\}$ and let
$\T_n=\{{\rm Range}(f):f\in\FF_\l(N,M)$ has domain $A_n\}$.
We will show that 
$\T=\bigcup\{\T_n:n\in\o\}$
is both well-founded and has depth $\al^*$.
As noted above, this is sufficient to conclude that
\hbox{$KC_\l(\KK)=\infty$}.
If we assume that $\T$  is well-founded then 
the family of maps 
$$f_\eta:A_{lg(\eta)}\ra M$$
for $\eta\in des(\al^*)$
defined by $f_\eta((\al,\eta_i))=(\al,\eta| i)$
witness that the depth of $\T$ is at least $\al^*$.

So it remains to show that $\T$ is well-founded.
The obvious distinction between the ordered multigraphs $I$ and $J$ is
that $J$ has an infinite, strictly increasing sequence $\lan\eta_n:n\in\o\ran$,
whereas $I$ does not.  
Suppose  that an elementary map 
$g:\bigcup\{A_n:n\in\o\}\ra M$ is given.
We will obtain  a contradiction by constructing an infinite
strictly increasing sequence in $des(\al^*)$.
The construction of this sequence proceeds in three stages.  First, since $\mu>|T_1|$ is regular,
for every $l\in\o$ there is an integer $n(l)$, an
$L_1$-term $\tao_l(x_1,\dots,x_{n(l)})$, a subset $X_l$ of $\mu$ of
size $\mu$, and functions $t_{l,m}:X_l\ra I$
such that for each $\b\in X_l$
$$g((\b,\eta_l))=\tao_l(d_l(\b)),$$
where $d_l(\b)=\lan t_{l,1}(\b),\dots,t_{l,n(l)}(\b)\ran$.
As notation, let $W=\{(l,m):l\in\o$ and $m\in[1,\dots,n(l)]\}$ and
for each $(l,m)\in W$, let $\al_{l,m}$ and $\eta_{l,m}$
be the functions with domain $X_l$ satisfying $$t_{l,m}(\b)=(\al_{l,m}(\b),\eta_{l,m}(\b)).$$
Next, we state two claims, whose proofs we defer until the end of the argument.

\medskip\par\noindent
{\bf Claim 1.}  There is a sequence $\lan Y_l:l\in\o\ran$ such that each $Y_l\sbq X_l$
has size $\mu$ and for each $k\in\o$
$$\hbox{qftp}(d_0(\b_0),\dots,d_{k-1}(\b_{k-1}))=\hbox{qftp}(d_0(\b_0'),\dots,d_{k-1}(\b_{k-1}'))$$
in the structure $(I,<,\lhd^*)$
for all sequences $\b_0<\dots<b_{k-1}$, $\b_0'<\dots<b_{k-1}'$ with $\b_l,\b'_l\in Y_l$ 
for each $l<k$.
\medskip\par\noindent
{\bf Claim 2.}  
For every $k>1$ there is a sequence $\lan m_l:l<k\ran$ and a
permutation $\sigma$ of $k$ such that
%for every $l<k$, $(l,m_l)$ is $\al$-trivial, and
$$t_{\sigma(0),m_{\sigma(0)}}(\b_{\sigma(0)})\, \lhd^* \,
t_{\sigma(1),m_{\sigma(1)}}(\b_{\sigma(1)})\, \lhd^*\dots\lhd^* \,
t_{\sigma(k-1),m_{\sigma(k-1)}}(\b_{\sigma(k-1)})$$
for every sequence $\b_0<\dots<\b_{k-1}$ with $\b_l\in Y_l$ for each $l<k$.

Given these two claims, it follows from K\"onig's Lemma (and the fact that the
permutation $\sigma$ is uniquely determined by the lengths of the $\eta^t$'s) that there
is an infinite sequence $\lan m_l:l\in\o\ran$ and a permutation $\sigma\in Sym(\o)$
such that, letting $\eta_l=\eta^{t_{\sigma(l),m_\sigma(l)}}$ for each $l\in\o$,
$$\eta_0(\b_{\sigma(0)})\lhd \eta_1(\b_{\sigma(1)})\lhd \dots$$
for all sequences $\b_0<\b_1<\dots$ satisfying $\b_l\in Y_l$ for each $l\in\o$.
But the existence of such a sequence is clearly impossible as each $\eta_l(\b)\in des(\al^*)$.
\endproof

Thus, to complete the proof of the theorem it suffices to prove the claims.
The proof of Claim~1 is tedious, but straightforward.
First, by trimming each of the sets $X_l$
we may assume that for each $(l,m)\in W$,
\begin{enumerate}
\item $\al_{l,m}$ is constant on $X_l$;
\item $\al_{l,m}(\b)=\b$ for all $\b\in  X_l$; or
\item $\{\al_{l,m}(\b):\b\in X_l\}$ is strictly increasing and disjoint from $X_l$.
\end{enumerate}
We call $(l,m)$ {\it $\al$-constant\/} if (1) holds 
%(we denote its value by $\al(l,m)$)
and  call $(l,m)$
{\it $\al$-trivial\/} if (2) holds.  
Similarly, we may assume that for
each $(l,m)\in W$,
\begin{itemize}
\item $lg(\eta_{l,m}(\b))$ is constant for all $\b\in X_l$ and
\item $\eta_{l,m}$ is constant on $X_l$ or  else $\{\eta_{l,m}(\b):\b\in X_l\}$
is strictly increasing (in lexicographic order).
%\item there is $\nu_{l,m}\in des(\al^*)$ such that $\nu_{l,m}\lhd \eta^{t^{\b,l}_m}$
%for all $\b\in X_l$ and:
%\begin{itemize}
%\item either $\nu_{l,m}= \eta^{t^{\b,l}_m}$ for all $\b\in X_l$
%or $\{\eta^{t^{\b,l}_m}(s):\b\in X_l\}$ is strictly increasing,
%where $s=lg(\nu_{l,m})$.
%\end{itemize}
\end{itemize}

Additionally, we may assume that for each pair $(l,m),(l,m')\in W$ with the same $l$,
the truth values of 
\begin{itemize}
\item ``$\al_{l,m}(\b) < \al_{l,m'}(\b)$'';
\item ``$\eta_{l,m}(\b) \lhd \eta_{l,m'}(\b)$''; 
\item ``$\eta_{l,m}(\b) <_{\rm lex} \eta_{l,m'}(\b)$''; and hence of
\item ``$t_{l,m}(\b) < t_{l,m'}(\b)$''
\end{itemize}
are  constant for all $\b\in X_l$.
By trimming each $X_l$ further, we may additionally assume
that for all pairs $m,m'\in[1,\dots n(l)]$,
the truth values of 
\begin{itemize}
\item ``$\al_{l,m}(\b_1) < \al_{l,m'}(\b_2)$'';
\item ``$\eta_{l,m}(\b_1) \lhd \eta_{l,m'}(\b_2)$'';
\item ``$\eta_{l,m}(\b_1) <_{\rm lex} \eta_{l,m'}(\b_2)$''; and hence of
\item ``$t_{l,m}(\b_1) < t_{l,m'}(\b_2)$'' 
\end{itemize}
are constant for all pairs $\b_1<\b_2$ from $X_l$.

So far, each of our trimmings has concentrated on a single set $X_l$.
However, to complete the proof of the claim, we must consider pairs of sets as well.
Fortunately, this presents no problem.  We illustrate one such reduction and leave the
other (virtually identical) reductions to the reader.
We claim that there are subsets $Y_l\sbq X_l$, each of size $\mu$, such that
for all $(l_1,m_1), (l_2,m_2)\in W$
the truth value of 
$$\hbox{``}\al_{l_1,m_1}(\b_1)< \al_{l_2,m_2}(\b_2)\hbox{ ''} \eqno(*)$$
is constant for all pairs $(\b_1,\b_2)$ satisfying 
$\b_1\in Y_{l_1}$, $\b_2\in Y_{l_2}$, and $\b_1<\b_2$.
To see this, let $C$ be the $\al$-constant pairs $(l,m)\in W$
and let $\d<\mu$ be the supremum of all 
$\al_{l,m}(\b)$ for $(l,m)\in C$, $\b\in X_l$.
By removing fewer than $\mu$ elements from each $X_l$, we may assume that
$\al_{l,m}(\b)>\d$ for all non-$\al$-constant $(l,m)\in W$ and all $\b\in X_l$.
It is now routine to inductively construct the sets $\{Y_l:l\in\o\}$ in $\mu$
steps so as to ensure
$$\al_{l_1,m_1}(\b_1)< \al_{l_2,m_2}(\b_2)$$
whenever $(l_1,m_1),(l_2,m_2)$ are not $\al$-constant, 
$l_1<l_2$, $\b_1\in Y_1$, $\b_2,\in Y_2$, and $\b_1<\b_2$.
Combining this with the earlier trimmings of the $X_l$'s establish ($*$).

Finally, we prove Claim 2.  This is the heart of the argument and
is where properties of the coloring $c$ are used.
Fix an integer $k>1$.  In light of 
Claim 1, it suffices to find a sequence $\lan m_l:l<k\ran$
and a permutation $\sigma$ of $k$ such that
$t_{\sigma(0),m_{\sigma(0)}}(\b_{\sigma(0)})\, 
\lhd^*\dots\lhd^* \,
t_{\sigma(k-1),m_{\sigma(k-1)}}(\b_{\sigma(k-1)})$
for {\it some\/} sequence 
$\b_0<\dots<\b_{k-1}$ with $\b_l\in Y_l$ for each $l<k$.
Consequently, we can trim the sets $Y_l$ still further.
As notation, let $W_k$ denote the {\it finite\/} set of all pairs
$(l,m)\in W$ with $l<k$.
For each $l<k$, let $h_l$ enumerate $Y_l$, i.e., 
$h_l(\d)=$ the $\d^{\rm th}$ element of $Y_l$.

By trimming each $Y_l$ for  $l<k$, we may additionally assume that:
\begin{itemize}
\item  The sets $Y_l$ are disjoint and $\al_{l,m}(\b)\not\in \bigcup_{l<k} Y_l$
unless $(l,m)$ is $\al$-trivial;
\item $\d_1<\d_2$ implies $h_l(\d_1)<h_{l'}(\d_2)$ for all $l,l'<k$;
\item For all pairs $(l,m),(l',m')\in W_k$ with $(l,m)$ $\al$-constant,
there is an integer $c^*(l,m,l',m')\in\o$ such that
$$c(\al_{l,m}(\b),\al_{l',m'}(\b'))=c^*(l,m,l',m')$$
for all distinct $\b,\b'$ from $Y_l,Y_l'$ respectively.
\end{itemize}
Let $C^*$ denote the (finite) set of all integers $c^*(l,m,l',m')$, where
the pairs $(l,m),(l',m')$ are from $W_k$ 
and $(l,m)$ is $\al$-constant.  Choose integers $p\in Z_k\setminus C^*$ and
$q\in Z_r$ for some $r>|W_k|$.  As notation, for each ordinal $\d\in\mu$, let
$$B_l(\d)=\{\al_{l,m}(h_l(\d)):(l,m)\in W,\ \hbox{$(l,m)$ not $\al$-constant}\}\cup\{h_l(\d)\}$$
For $\delbar=\d_0<\d_1<\dots<\d_{k-1}$, let $B(\delbar)=\bigcup_{l<k} B_l(\d_l)$.
By trimming the sets $Y_l$, $l<k$ still further, we may assume
that the order type of $B(\delbar)$ is constant
among all increasing $k$-tuples $\delbar$.
Thus, by employing Lemma~\ref{useful}, we can choose two increasing $k$-tuples 
$\delbar^0$ and $\delbar^1$ satisfying:
\begin{itemize}
\item $c(\al,\b)=q$ for all $\al,\b\in B(\delbar^0)$; and
\item $c(\al,\b)=q$ for all $\al,\b\in B(\delbar^1)$
EXCEPT that \hbox{$c(h_i(\d_i^1),h_j(\d_j^1))=p$} for all $i\neq j$.

\end{itemize}

As notation, let $\nu_l=h_l(\d^0_l)$, 
$\nubar=\nu_0<\dots<\nu_{k-1}$,
and $D(\nubar)=\{t_{l,m}(\nu_l):(l,m)\in W_k\}$.  Dually, let
$\b_l=h_l(\d^1_l)$, 
$\betabar=\beta_0<\dots<\beta_{k-1}$, and
$D(\betabar)=\{t_{l,m}(\b_l):(l,m)\in W_k\}$.  

Now, working in the multigraph $J$, 
%as $\rho_0\lhd\rho_1 \lhd\dots\lhd\rho_{k-1}$,
$$J\models \neg R_k((\nu_0,\rho_0),\dots,(\nu_{k-1},\rho_{k-1}))\wedge 
R_k((\b_0,\rho_0),\dots,(\b_{k-1},\rho_{k-1})),$$
so 
$$N\models \neg\phi_k( (\nu_0,\rho_0),\dots,(\nu_{k-1},\rho_{k-1}))\wedge
\phi_k( (\b_0,\rho_0),\dots,(\b_{k-1},\rho_{k-1})).$$
Hence, by the elementarity of the map $g$,
$$ M\models  \neg \phi_k(\tao_0(d_0(\nu_0)),\dots,\tao_{k-1}(d_{k-1}(\nu_{k-1})))\ \hbox{and}$$
$$ M\models  \phi_k(\tao_0(d_0(\b_0)),\dots,\tao_{k-1}(d_{k-1}(\b_{k-1}))).$$

It follows from Proposition~\ref{coding}
that the discrepancy in $\phi_k$
implies that ${\rm qftp}(D(\nubar))\neq {\rm qftp}(D(\betabar))$
in the ordered multigraph $(I,<,R_n)_{n\in\o}$.
However, since 
${\rm qftp}(D(\nubar))= {\rm qftp}(D(\betabar))$
in the structure $(I,<,\lhd^*)$, the sets $D(\nubar)$ and $D(\betabar)$
must differ on some $R_n$.  This difference
can only be explained by a discrepancy of the function $c$
on some pairs of elements from the sets $B(\nubar)$ and $B(\betabar)$.
Since $c$ can only attain the values of $p$ and $q$ on 
pairs from $B(\nubar)$ and
$B(\betabar)$, our choice of $p$ and $q$
implies that $R_k$ is the only relation that can differ
between $D(\nubar)$ and $D(\betabar)$.
Thus, there are sequences $s_0,\dots,s_{k-1}\in D(\betabar)$ and $s_0',\dots,s_{k-1}'
\in D(\nubar)$ of corresponding elements such that
$$I\models R_k(s_0,\dots,s_{k-1})\wedge\neg R_k(s_0',\dots,s_{k-1}').$$
In particular, $s_0\lhd^*\dots\lhd^* s_{k-1}$ and for all $i<j<k$ we have
$\al^{s_i}\neq\al^{s_j}$ and $c(\al^{s_i},\al^{s_j})=p$.
As each $s_i\in D(\betabar)$, there are functions $l,m$ with domain $k$ such 
that $$s_i=t_{l(i),m(i)}(\b_{l(i)}).$$
Now fix $i<k$.  Since $k>1$ and $c(\al^{s_i},\al^{s_j})=p\not\in C^*$ for any $j\neq i$,
the pair $(l(i),m(i))$ is not $\al$-constant.  As well, the choice of the coloring
of $B(\betabar)$ ensures that $\al^{s_i}=h_r(\d^1_r)\in Y_r$ for some $r<k$.
Thus, the disjointness of the $Y_l$'s imply that $r=l(i)$ and the pair $(l(i),m(i))$ is
$\al$-trivial.  That is, $\al^{s_i}=\b_{l(i)}$.
Further, since $\al^{s_i}\neq\al^{s_j}$ whenever $i\neq j$, 
the function $l$ must be a permutation of the set $k$.
So, letting  $m_i=m(l^{-1}(i))$ and $\sigma=l$, the sequence $\lan m_i:i<k\ran$ and
permutation $\sigma$ are as desired.  \endproof

%
%
%\bibliography{mclbis}
%\bibliographystyle{plain}

\end{document}